\documentclass[11pt]{article}
\usepackage{epsfig}

\def\be{\begin{equation}}
\def\ee{\end{equation}}
\def\bea{\begin{eqnarray}}
\def\eea{\end{eqnarray}}
\def\bes{\begin{eqnarray*}}
\def\ees{\end{eqnarray*}}

\def\nn{\nonumber}
\def\<{\langle}
\def\>{\rangle}
\def\lb{\label}

\def\R{{\bf R}}
\def\C{{\bf C}}
\def\Z{{\bf Z}}

\def\N{{\bf N}}
\def\U{{\bf U}}

\def\Q{{\bf Q}}

\def\aa{{\alpha}}

\def\ga{{\gamma}}

\def\th{{\theta}}

\def\Lm{{\Lambda}}

\def\sg{{\sigma}}

\def\rank{{\rm rank}}
\def\Sp{{\rm Sp}}

\def\ol{\overline}

\def\ol#1{\overline{#1}}  

\def\hb{\vrule height0.18cm width0.14cm $\,$}

\def\ol#1{\overline{#1}}  

\title{Multiplicity and stability of closed geodesics on bumpy
Finsler $3$-spheres}

\author{Huagui Duan$^{1}$,\thanks{Partially supported by NNSF and RFDP
of MOE of China. E-mail: duanhuagui@163.com}
\qquad Yiming Long$^{1,2}$\thanks{Partially supported by the 973 Program of
MOST, Yangzi River Professorship, NNSF, MCME, RFDP, LPMC of MOE of
China, S. S. Chern Foundation, and Nankai University.
E-mail: longym@nankai.edu.cn }\\ \\
$^{1}$ Chern Institute of Mathematics\\
$^{2}$ Key Lab of Pure Mathematics and Combinatorics of Ministry of
Education\\ Nankai University, Tianjin 300071\\ The People's
Republic of China\\ }
\date{}

\begin{document}

\maketitle

\begin{abstract}
{\it  We prove that for every $\Q$-homological Finsler 3-sphere
$(M,F)$ with a bumpy and irreversible metric $F$, either there
exist two non-hyperbolic prime closed geodesics, or there exist at
least three prime closed geodesics.}
\end{abstract}

\renewcommand{\theequation}{\thesection.\arabic{equation}}
\renewcommand{\thefigure}{\thesection.\arabic{figure}}

\setcounter{equation}{0}
\section{Introduction and the main result}

Let $M$ be a smooth manifold with a Finsler metric $F$. A continuous
curve $c:[0,1]\to M$ on a Finsler manifold $(M,F)$ is a closed
geodesic, if $c$ is closed and is the shortest curve connecting any
two points on the image of $c$ which are close enough (cf.
\cite{BCS1} and \cite{She1}). As usual on any Finsler manifold
$M=(M,F)$, a closed geodesic $c:S^1=\R/\Z\to M$ is {\it prime}, if
it is not a multiple covering (i.e., iteration) of any other closed
geodesics. Here the $m$-th iteration $c^m$ of $c$ is defined by
$c^m(t)=c(mt)$ for $m\in\N$. The inverse curve $c^{-1}$ of $c$ is
defined by $c^{-1}(t)=c(1-t)$ for $t\in \R$. We call two prime
closed geodesics $c$ and $d$ {\it distinct} if there is no
$\th\in(0,1)$ such that $c(t)=d(t+\th)$. We shall omit the word
``distinct" for short when we talk about more than one prime closed
geodesics. A closed geodesic $c$ is {\it non-degenerate} if $1$ is
not an eigenvalue of the linearized Poincar\'e map $P_c$ of $c$. $c$
is {\it hyperbolic} if $\sg(P_c)\cap \U=\emptyset$, or {\it
elliptic} if $\sg(P_c)\subset\U$, where we denote by
$\U=\{z\in\C\,|\,|z|=1\}$. A Finsler metric $F$ is {\it bumpy} if
all closed geodesics and their iterations are all non-degenerate.

Note that by the classical theorem of Lyusternik-Fet \cite{LyF1},
there exists at least one closed geodesic on every compact
Riemannian as well as Finsler manifold. In \cite{Kat1} of 1973, Katok
constructed his famous irreversible Finsler metrics on $S^{2n}$
respectively $S^{2n-1}$ with precisely $2n$ prime closed geodesics all
of which are elliptic (cf. also \cite{Zil1}). Based on this result in
\cite{Ano1} of 1974, Anosov conjectured that the lower bound of the number
of prime closed geodesics on $S^n$ is $2[\frac{n}{2}]$.

We are only aware of a few results on the multiplicity and stability
of closed geodesics on Finsler spheres. In 1965, Fet in
\cite{Fet1} proved that there exist at least two distinct closed
geodesics on every reversible bumpy Finsler manifold $(M,F)$. In 1989,
Rademacher in \cite{Rad1} proved that there exist at least
two elliptic closed geodesics on every bumpy Finsler 2-sphere with
finitely many prime closed geodesics. In 2003, Hofer, Wysocki and Zehnder
in \cite{HWZ1} proved that there exist either two or infinitely many
prime closed geodesics on every bumpy Finsler 2-sphere if the stable
and unstable manifolds of every hyperbolic closed geodesic intersect
transversally. In \cite{Rad3} of 2005, Rademacher obtained existences and
stability of closed geodesics on Finsler $S^n$ under pinching
conditions which generalizes results in \cite{BTZ1} and
\cite{BTZ2} of Ballmann, Thorbergsson and Ziller in
1982-83 on Riemannian manifolds. Recently, Bangert and Long
in \cite{BaL1} proved that there exist always at least two prime
closed geodesics on every Finsler 2-sphere $(S^2,F)$ which answers
Anosov's conjecture for $S^2$. More recently Duan and Long in
\cite{DuL1} and Rademacher in \cite{Rad4} proved independently
that there exist at least two distinct prime closed geodesics on
every bumpy Finsler $n$-sphere $(S^n,F)$ with $n\ge 3$.

Note that in \cite{Hin1} of 1984, Hingston's result specially
implies that on every Riemannian sphere $S^n$, if all the closed
geodesics are hyperbolic, then there exist infinitely many
geometrically distinct closed geodesics. In \cite{Rad1} of 1989,
Rademacher proved that on every even dimensional bumpy Finsler sphere
$S^{2n}$, if there are only finitely many prime closed geodesics, then
at least one of them is non-hyperbolic. In the recent \cite{LoW1}, Long
and Wang proved that if there exist precisely two prime closed geodesics
on a Finsler $S^2$, both of them must be irrationally elliptic.

Note that Long in \cite{Lon3} conjectured that the number of prime
closed geodesics for $(S^3,F)$ may belong to
$\{2,3,4\}\bigcup\{+\infty\}$. This paper is devoted to the proof of
the following main result which is related to this conjecture.

{\bf Theorem 1.1.} {\it For every $\Q$-homological Finsler
$3$-sphere $(M,F)$ with a bumpy and irreversible metric $F$, either
there exist precisely two non-hyperbolic prime closed geodesics, or
there exist at least three distinct prime closed geodesics. }

Theorem 1.1 is a consequence of the slightly stronger version of
Theorem 3.6 below. Our proof of these theorems relies mainly on the
following ingredients: The precise index iteration formulae of Long,
the common index jump theorem of Long and Zhu, the Morse
inequalities, and Rademacher's mean index identity for closed
geodesics, together with some new techniques relating local and
global information. Because the proof for $\Q$-homological
$3$-spheres is the same as that for $S^3$, we carry out below the
proof only for $(S^3,F)$ with a bumpy and irreversible Finsler
metric $F$. The main idea is that assuming the existence of
precisely two prime closed geodesics on $(S^3,F)$ and at least one
of them being hyperbolic, we shall derive a contradiction via the
above mentioned tools.

In this paper, let $\N$, $\N_0$, $\Z$, $\Q$, $\R$, and $\C$ denote
the sets of positive integers, non-negative integers, rational
numbers, real numbers, and complex numbers respectively. We denote
by $[a]=\max\{k\in\Z\,|\,k\le a\}$ for any $a\in\R$. When $S^1$
acts on a topological space $X$, we denote by $\overline{X}$ the
quotient space $X/S^1$. We use only singular homology modules with
$\Q$-coefficients. For terminologies in algebraic topology we
refer to \cite{GrH1}.

\setcounter{equation}{0}
\section{Preliminary results on closed geodesics}

\subsection{Critical modules of iterations of closed geodesics}

Let $M=(M,F)$ be a compact Finsler manifold $(M,F)$, the space
$\Lambda=\Lambda M$ of $H^1$-maps $\gamma:S^1\rightarrow M$ has a
natural structure of Riemannian Hilbert manifolds on which the
group $S^1=\R/\Z$ acts continuously by isometries, cf.
\cite{Kli2}, Chapters 1 and 2. This action is defined by
$(s\cdot\gamma)(t)=\gamma(t+s)$ for all $\gamma\in\Lm$ and $s,
t\in S^1$. For any $\gamma\in\Lambda$, the energy functional is
defined by
\be E(\gamma)=\frac{1}{2}\int_{S^1}F(\gamma(t),\dot{\gamma}(t))^2dt.
\lb{2.1}\ee
It is $C^{1,1}$ (cf. \cite{Mer1}) and invariant under the $S^1$-action. The
critical points of $E$ of positive energies are precisely the closed geodesics
$\gamma:S^1\to M$. The index form of the functional $E$ is well
defined along any closed geodesic $c$ on $M$, which we denote by
$E''(c)$ (cf. \cite{She1}). As usual, we denote by $i(c)$ and
$\nu(c)$ the Morse index and nullity of $E$ at $c$. In the
following, we denote by
\be \Lm^\kappa=\{d\in \Lm\;|\;E(d)\le\kappa\},\quad \Lm^{\kappa-}=\{d\in \Lm\;|\; E(d)<\kappa\},
  \quad \forall \kappa\ge 0. \nn\ee
For a closed geodesic $c$ we set $ \Lm(c)=\{\ga\in\Lm\mid E(\ga)<E(c)\}$.

For $m\in\N$ we denote the $m$-fold iteration map
$\phi_m:\Lambda\rightarrow\Lambda$ by $\phi_m(\ga)(t)=\ga(mt)$, for all
$\,\ga\in\Lm, t\in S^1$, as well as $\ga^m=\phi_m(\gamma)$. If $\gamma\in\Lambda$
is not constant then the multiplicity $m(\gamma)$ of $\gamma$ is the order of the
isotropy group $\{s\in S^1\mid s\cdot\gamma=\gamma\}$. For a closed geodesic $c$,
the mean index $\hat{i}(c)$ is defined as usual by
$\hat{i}(c)=\lim_{m\to\infty}i(c^m)/m$. Using singular homology with rational
coefficients we consider the following critical $\Q$-module of a closed geodesic
$c\in\Lambda$:
\be \overline{C}_*(E,c)
   = H_*\left((\Lm(c)\cup S^1\cdot c)/S^1,\Lm(c)/S^1\right). \lb{2.3}\ee

The following results of Rademacher will be used in our proofs below.

{\bf Proposition 2.1.} (cf. Satz 6.11 of \cite{Rad2} ) {\it Let $c$ be a
prime closed geodesic on a bumpy Finsler manifold $(M,F)$. Then there holds}
$$ \overline{C}_q( E,c^m) = \left\{\matrix{
     \Q, &\quad {\it if}\;\; i(c^m)-i(c)\in 2\Z\;\;{\it and}\;\;
                   q=i(c^m),\;  \cr
     0, &\quad {\it otherwise}. \cr}\right.  $$

{\bf Definition 2.2.} (cf. Definition 1.6 of \cite{Rad1}) {\it For a
closed geodesic $c$, let $\ga_c\in\{\pm\frac{1}{2},\pm1\}$ be the
invariant defined by $\ga_c>0$ if and only if $i(c)$ is even, and
$|\ga_c|=1$ if and only if $i(c^2)-i(c)$ is even. }

{\bf Proposition 2.3.} (cf. Theorem 3.1 of \cite{Rad1}) {\it Let
$c_k, k=1,2,\cdots,r$ prime closed geodesics on a  bumpy Finsler
$3$-sphere. Then the average indices $\hat{i}(c_k)$ and the
invariants $\ga_{c_k}$ satisfy the identity }
\be \sum_{k=1}^r\frac{\ga_{c_k}}{\hat{i}(c_k)}=1.    \lb{2.4}\ee

Let $(X,Y)$ be a space pair such that the Betti numbers
$b_i=b_i(X,Y)=\dim H_i(X,Y;\Q)$ are finite for all $i\in \Z$. As
usual the {\it Poincar\'e series} of $(X,Y)$ is defined by the
formal power series $P(X, Y)=\sum_{i=0}^{\infty}b_it^i$. We need
the following well known results on Betti numbers and the Morse
inequality for $\overline{\Lm}\equiv
\overline{\Lm} S^3$ and $\ol{\Lm}^0=\ol{\Lambda}^0S^3
=\{{\rm constant\;point\;curves\;in\;}S^3\}\cong S^3$.

{\bf Proposition 2.4.} (cf. Remark 2.5 of \cite{Rad1}) {\it The
Poincar\'e series is given by \bea P(\ol{\Lm}S^3,\ol{\Lm}^0S^3)(t)
&=&t^2\left(\frac{1}{1-t^2}+\frac{t^2}{1- t^2}\right)  \nn\\
&=& t^2(1+t^2)(1+t^2+t^4+\cdots) = t^2+2t^4+2t^6+\cdots, \nn\eea
which yields}
\be {b}_q = {b}_q(\ol{\Lm}S^3,\ol{\Lm}^0 S^3)\;
 = \rank H_q(\ol{\Lm} S^3,\ol{\Lm}^0 S^3 )
 = \;\;\left\{\matrix{
    1,&\quad {\it if}\quad q=2,  \cr
    2,&\quad {\it if}\quad q=2k+2,\quad k\in \N,  \cr
    0 &\quad {\it otherwise}. \cr}\right. \lb{2.5}\ee

{\bf Proposition 2.5.} (cf. Theorem I.4.3 of \cite{Cha1}, Theorem
6.1 of \cite{Rad2}) {\it Suppose that there exist only finitely
many prime closed geodesics $\{c_j\}_{1\le j\le k}$ on a Finsler
$3$-sphere $(S^3, F)$. Set
$$ M_q =\sum_{1\le j\le k,\; m\ge 1}\dim{\ol{C}}_q(E, c^m_j), \quad \forall q\in\Z. $$
Then for every integer $q\ge 0$ there holds }
\bea
M_q - M_{q-1} + \cdots +(-1)^{q}M_0
&\ge& b_q - b_{q-1}+ \cdots + (-1)^{q}b_0, \lb{2.6}\\
M_q &\ge& b_q. \lb{2.7}\eea

\subsection{Classification of non-degenerate closed geodesics on $S^3$ }

We introduce some notations in \cite{Lon2} here. Given any
two real matrices of the square block form
$$ M_1=\left(\matrix{A_1 & B_1\cr C_1 & D_1\cr}\right)_{2i\times 2i},
  \qquad M_2=\left(\matrix{A_2 & B_2\cr C_2 & D_2\cr}\right)_{2j\times 2j},$$
we define the $\diamond$-sum of $M_1$ and $M_2$ to be the $2(i+j)\times2(i+j)$
matrix $M_1\diamond M_2$ given by
$$ M_1\diamond M_2=\left(\matrix{A_1 & 0 & B_1 & 0 \cr
                                   0 & A_2 & 0& B_2\cr
                                 C_1 & 0 & D_1 & 0 \cr
                                   0 & C_2 & 0 & D_2}\right),$$
and $M_1^{\diamond k}$ to be the $k$-times $\diamond$-sum of $M_1$. In the
following, let
\bes N({\aa},{B})&=&\left(\matrix{\cos{\aa}
                     & -\sin{\aa} & {b}_{1} & {b}_{2} \cr \sin{\aa} & \cos{\aa} &
{b}_{3} & {b}_{4}\cr 0 & 0 & \cos{\aa} & -\sin{\aa} \cr 0 & 0 &
\sin{\aa} & \cos{\aa}
                  }\right),\\[6pt] R(\th)&=&\left(\matrix{\cos\th &
-\sin\th\cr \sin\th & \cos\th\cr}\right)\quad\mbox{and}\quad
H(d)=\left(\matrix{ d & 0\cr 0 & 1/d\cr}\right),\ees
where
$B=\left(\matrix{b_{1} & b_{2}\cr b_{3} & b_{4}\cr}\right)$ with
$({b}_{1},{b}_{2},{b}_{3},{b}_{4})\in\R^4$, ${\th}/\pi$, and
$\aa/\pi\in (0, 2)\setminus\Q$, and $d>0$ or $d<0$.

In \cite{Lon1} of 2000, Long defined the {\it homotopy set}
$\Omega(M)$ and the {\it homotopy component} $\Omega^0(M)$ of $M$ in
the symplectic group $\Sp(2n)$ by
\bes\Omega(M)=\{N\in \Sp(2n)\ |\ \sigma(N)\cap\U=\sigma(M)\cap\U\equiv
\Gamma\;\mbox{and}\;\nu_{\omega}(N)=\nu_{\omega}(M)\ \forall
\omega\in\Gamma\},\ees
where $\sigma(M)$ denotes the spectrum of $M$ and
$\nu_{\omega}(M)\equiv\dim_{\C}\ker_{\C}(M-\omega I)$ for all $\omega\in\U$.
Then $\Omega^0(M)$ is defined to be the path connected component of
$\Omega(M)$ containing $M$ (cf. also Section 1.8 of \cite{Lon2}).

Let $c$ be a closed geodesic on a Finsler sphere $(S^3,F)$.
Denote the linearized Poincar\'e map of $c$ by $P_c$. Note that the
index iteration formulae in \cite{Lon1} (cf. also \cite{Lon2}) work for
Morse indices of iterated closed geodesics on Finsler manifolds (cf.
\cite{LLo1}). Suppose that all iterations $c^m$ of $c$ are non-degenerate.
Then by Theorems 8.1.4 to 8.1.7 and Theorem 8.3.1 of \cite{Lon2}, we
have the following classification of non-degenerate closed
geodesics, i.e., there exists a path $f_c\in C([0,1],\Omega^0(P_c))$
such that $f_c(0)=P_c$ and $f_c(1)$ have the following forms.

{\bf NCG-1.} {\it $f_c(1)=R(\th_1)\diamond R(\th_2)$.}

In this case, by Theorems 8.1.7 and 8.3.1 of \cite{Lon2}, we have $i(c)=2p$
for some $p\in\N_0$, and
\be i(c^m) = 2m(p-1) +2\sum_{i=1}^{2}\left[\frac{m\th_i}{2\pi}\right]+2 , \quad
       \nu(c^m)=0, \qquad \forall \ m\ge 1.\lb{2.8} \ee

{\bf NCG-2.} {\it $f_c(1)=R(\th)\diamond H(d)$.}

In this case, by Theorems 8.1.6, 8.1.7 and 8.3.1 of \cite{Lon2}, we have
$i(c)=p$ for some $p\in\N_0$, and
\be i(c^m) = m(p-1) +2\left[\frac{m\th}{2\pi}\right]+1 , \quad
    \nu(c^m)=0, \qquad \forall \ m\ge 1.\lb{2.9} \ee

{\bf NCG-3.} {\it $f_c(1)=  H(d_{1})\diamond H(d_2)$.}

In this case, by Theorems 8.1.6 and 8.3.1 of \cite{Lon2}, we have
$i(c)=p$ for some $p\in\N_0$, and
\be i(c^m) = mp , \quad \nu(c^m)=0, \qquad \forall \ m\ge 1.\lb{2.10} \ee

{\bf NCG-4.} {\it $f_c(1)=N(\aa,B)$.}

In this case, by Theorems 8.2.3, 8.2.4 and 8.3.1 of \cite{Lon2}, we have
$i(c)=2p$ for some $p\in\N_0$, and
\be i(c^m) = 2mp, \quad \nu(c^m)=0, \qquad \forall \ m\ge 1.\lb{2.11} \ee

\setcounter{equation}{0}
\section{Proof of main theorems}

Firstly, we list below three auxiliary results. By Theorem 1.2 in \cite{DuL1} (cf.
also \cite{Rad4}), there exist at least two prime closed geodesics on every bumpy
Finsler sphere $(S^n,F)$ for $n\ge 2$. As also noticed in Rademacher's preprint
\cite{Rad4}, we point out that the {\it common index jump theorem} due to Long
and Zhu in \cite{LoZ1} also works for closed geodesics on Finsler manifolds
(cf. Remark 12.2.5 of \cite{Lon2}). In the bumpy case, we have the following
consequence.

{\bf Theorem 3.1} (cf. Theorem 4.3 of \cite{LoZ1} and Theorem 11.2.1
of \cite{Lon2}) {\it Let $c$ be a closed geodesic on a compact bumpy Finsler manifold
$(M,F)$ with $\hat{i}(c)>0$. Then there exist infinitely many $k\in\N$ such that }
\be i(c^{2k+1})-i(c^{2k-1})=2i(c). \lb{3.1}\ee

{\bf Proof.} For readers conveniences, we sketch the proof here. By Theorem 4.3 of
\cite{LoZ1} (cf. (11.2.4) and (11.2.5) in Theorem 11.2.1 of
\cite{Lon2}) we obtain for infinitely many $k\in\N$,
$$ i(c^{2k+1})-i(c^{2k-1}) - 2i(c) = 2S_{P_c}^+(1) + \nu(c^{2k-1})-\nu(c), $$
where $S_{P_c}^+(1)$ is the positive splitting number of $P_c$ at $1$. Because
$(M,F)$ is bumpy, all terms in the right hand side of the above identity are zero.
\hfill\hb

{\bf Lemma 3.2.} {\it Let $(M,F)$ a bumpy Finsler manifold with only finitely many
prime closed geodesics. If $c$ is a closed geodesic on $(M,F)$ which is not an
absolute minimum of the energy functional $E$ in its free homotopy class, the mean
index of the closed geodesic $c$ must satisfy $\hat{i}(c)> 0$. This holds
always when $M$ is simply connected, specially a sphere. }

{\bf Proof.} It is well known that the Morse index sequence $i(c^m)$
either tends to $+\infty$ asymptotically linearly or $i(c^m)=0$ for
all $m\ge 1$. Therefore $\hat{i}(c)=0$ if and only if $i(c^m)=0$ for
all $m\ge 1$. A crucial point in the proof of this fact is the
Property (2) in Proposition 1.3 in \cite{Bot1} of Bott (cf. also
Lemma 1 of \cite{GrM1}).

Now let $c$ be a prime closed geodesic on $(M,F)$ with $\hat{i}(c)=0$ which is not
an absolute minimum of $E$. Because $(M,F)$ is bumpy, every iteration $c^m$ of $c$
is homologically visible, by Theorem 3 of \cite{BaK1} there must exist infinitely
many prime closed geodesics on $(M,F)$. \hfill\hb

{\bf Lemma 3.3} {\it Let $(M,F)$ a Finsler manifold with only finitely many
prime closed geodesics. If the Morse type numbers $M_{2k-1}=0$ for all $k\in\N$,
then $M_q=b_q$ holds for all $q\in\N_0$. Specially, on a bumpy Finsler $(S^3,F)$
with finitely many prime closed geodesics, if the Morse indices of iterations
of these prime closed geodesics are all even, $M_q=b_q$ holds for all $q\in\N_0$.}

{\bf Proof.} It follows directly from Proposition 2.1 and the Morse
inequality. \hfill\hb

\smallskip

In this section we denote the contribution of iterations of prime
closed geodesics $c_i$ to the Morse type number $M_q$ by $M_q(i)$
for $i=1, 2$ and $q\ge 0$ below. The following lemma is crucial in
the proof of our main theorems.

\smallskip

{\bf Lemma 3.4.} {\it Let $S^3=(S^3,F)$ be a bumpy Finsler sphere with
precisely two prime closed geodesics $c_1$ and $c_2$. Suppose that $c_1$
and $c_2$ do not belong to the classes $\{\mbox{NCG-1},\mbox{NCG-2}\}$
simultaneously. Then at least one of the two closed geodesics must
satisfy $i(c)=2$. }

{\bf Proof.} By the Morse inequality and Proposition 2.4, we have
$M_2\ge b_2=1$. Thus at least one of the two closed geodesics must
have Morse index $i(c)\le 2$. Without loss of generality, let
$i(c_1)\le 2$. Assume the Lemma 3.4 does not hold. Then we have
$0\le i(c_1)\le 1$ and $i(c_2)\neq 2$. If $c_i$ for $i=1$ or $2$ belongs to
the class NCG-3 or NCG-4, then $i(c_i)>0$ by Lemma 3.2. So by
(\ref{2.11}), $c_1$ does not belong to NCG-4. Next we carry out our
proof in four cases according to the value of $i(c_1)$ and the
classification of $c_1$.

{\bf Case} I: {\it $c_1\in$NCG-3 with $i(c_1)=1$.}

By Proposition 2.1, there holds
$M_{2k+1}(1)=1,\,M_{2k}(1)=0,k\in\N_0$. We continue our study in 5
subcases (i)-(v):

(i) If $c_2\in\{\mbox{NCG-3},\;\mbox{NCG-4}\}$ with $i(c_2)=2p\ge4$.
Then by Definition 2.2, we have $\ga_{c_1}=-\frac{1}{2}$ and
$\ga_{c_2}=1$. Hence by Proposition 2.3 we obtain
$0>\frac{1}{2p}-\frac{1}{2}=\frac{\ga_{c_2}}{\hat{i}(c_2)}-\frac{\ga_{c_1}}{\hat{i}(c_1)}=1$
contradiction!

(ii) If $c_2\in$NCG-3 with $i(c_2)$ odd, then by Proposition 2.1,
Morse type numbers $M_{2k}(2)=0,k\in\N_0$. So $M_2=0$, which
contradicts to $M_2\ge b_2=1$ by Propositions 2.4 and 2.5.

(iii) If $c_2\in$NCG-2 with $i(c_2)=p\ge 0$, then
$\hat{i}(c_2)=p-1+\frac{\th}{\pi}$ is an irrational number by the
definition of $\th$. Hence
$\sum_{i=1}^2\frac{\ga_{c_i}}{\hat{i}(c_i)}$ is also an irrational
number. However, by Proposition 2.3 we have
$\sum_{i=1}^2\frac{\ga_{c_i}}{\hat{i}(c_i)}=1$ contradiction!

(iv) If $c_2\in$NCG-1 with $i(c_2)=2p\ge4$, then by Proposition 2.1,
$M_2(2)=0$. So $M_2=M_2(1)+M_2(2)=0$, which contradicts to $M_2\ge
b_2=1$ by Propositions 2.4 and 2.5.

(v) If $c_2\in$ NCG-1 with $i(c_2)=0$, then we have
$\hat{i}(c_2)=\frac{\th_1}{\pi}+\frac{\th_2}{\pi}-2,\;\hat{i}(c_1)=1$
by (\ref{2.8}) and (\ref{2.10}) in Section 2. By Definition 2.2 we
obtain $\ga_{c_1}=-\frac{1}{2}$ and $\ga_{c_2}=1$. Then by
Proposition 2.3, we obtain the identity \be
\frac{\th_1}{\pi}+\frac{\th_2}{\pi}=\frac{8}{3}. \lb{3.2}\ee On the
other hand, in this subcase we have $M_{2k-1}(2)=0$ for $k\in\N$.
Thus \be M_{2k-1}=M_{2k-1}(1)+M_{2k-1}(2)=1,\quad \forall k\in\N.
\lb{3.3}\ee

\smallskip

{\bf Claim 1}: {\it $i(c_2)=0$, $i(c^2_2)=i(c_2^3)=2$,
$i(c_2^4)=i(c_2^5)= i(c_2^6)=4$ and $i(c_2^7)=6$.}

\smallskip

In fact, by (\ref{3.2}) we have $\sum_{i=1}^2\frac{\th_i}{\pi}<3$.
Noting that $i(c_2^m)$ is even, it follows from
$i(c_2^2)=2\sum_{i=1}^{2}\left[\frac{\th_i}{\pi}\right]-2\le2$ that
$i(c_2^2)=0$ or $2$. If $i(c_2^2)=0$, then $M_0=M_0(2)\ge2$. However
it follows from Propositions 2.4, 2.5 and (\ref{3.3}) that $-1\ge
M_1-M_0\ge b_1-b_0=0$, which is a contradiction. So $i(c_2^2)=2$. By
(\ref{3.2}), we have $\frac{3\th_1}{2\pi}+\frac{3\th_2}{2\pi}=4$.
Since both $\frac{3\th_1}{2\pi}$ and $\frac{3\th_2}{2\pi}$ are
irrational numbers by the definition of $\th_i,i=1,2$, we have
$\sum_{i=1}^2[\frac{3\th_i}{2\pi}]\le3$. And so
$i(c_2^3)=2\sum_{i=1}^{2}\left[\frac{3\th_i}{2\pi}\right]-4\le2$. By
the similar argument it yields $i(c_2^3)=2$.

By (\ref{3.2}), we have
$\sum_{i=1}^2\frac{k\th_i}{2\pi}=\frac{4k}{3}=k+1+\frac{k-3}{3}$ for
$k=4,5,6$. Since both $\frac{k\th_1}{2\pi}$ and
$\frac{k\th_2}{2\pi}$ are irrational numbers, we have
$\sum_{i=1}^2[\frac{k\th_i}{2\pi}]\le k+1$. And so
$i(c_2^k)=2\sum_{i=1}^{2}\left[\frac{k\th_i}{2\pi}\right]-2k+2\in\{0,2,4\}$.
If $i(c_2^k)=0$ or $2$ for some $k\in\{4,5,6\}$, then $M_0\ge2$ or
$M_2\ge3$. By Propositions 2.4, 2.5 and (\ref{3.3}) we have
$-2\ge M_3-M_2+M_1-M_0\ge b_3-b_2+b_1-b_0=-1$, which is a contradiction. So
$i(c_2^4)=i(c_2^5)= i(c_2^6)=4$.

By (\ref{3.2}), we have
$\frac{7\th_1}{2\pi}+\frac{7\th_2}{2\pi}=\frac{28}{3}$. Since both
$\frac{7\th_1}{2\pi}$ and $\frac{7\th_2}{2\pi}$ are irrational
numbers, we have $\sum_{i=1}^2[\frac{7\th_i}{2\pi}]\le9$. And so
$i(c_2^7)=2\sum_{i=1}^{2}\left[\frac{7\th_i}{2\pi}\right]-12\le 6$.
By the above argument we then get $i(c_2^7)=6$. Claim 1 is proved.

\smallskip

Next we will estimate values of
$\{\frac{\th_1}{\pi},\frac{\th_2}{\pi}\}$ by analyzing
$i(c_2^m)$ for $m=2,4,5,6,7$ respectively.

By Claim 1, we obtain
$2=i(c_2^2)=2\sum_{i=1}^{2}\left[\frac{\th_i}{\pi}\right]-2$ which
implies
\be \sum_{i=1}^{2}\left[\frac{\th_i}{\pi}\right]=2.\lb{3.4}\ee
By (\ref{3.2}), without loss of generality, we get
$\frac{\th_1}{\pi}<\frac{4}{3}<\frac{\th_2}{\pi}$. By the definition
of $\th_2$, we have $\frac{\th_2}{\pi}<2$. Then by (\ref{3.4}),
$\frac{\th_1}{\pi}>1$. In summary, we obtain
\be 1<\frac{\th_1}{\pi}<\frac{4}{3}\quad\mbox{and}
      \quad \frac{4}{3}<\frac{\th_2}{\pi}<2. \lb{3.5}\ee

By Claim 1 we have
$4=i(c_2^4)=2\sum_{i=1}^{2}\left[2\frac{\th_i}{\pi}\right]-6$, which
implies
\be \sum_{i=1}^{2}\left[2\frac{\th_i}{\pi}\right]=5. \lb{3.6}\ee
Hence by (\ref{3.5}) we have
\be 2<2\frac{\th_1}{\pi}<\frac{8}{3}\quad\mbox{and}\quad
             \frac{8}{3}<2\frac{\th_2}{\pi}<4.  \lb{3.7}\ee
By (\ref{3.7}) we have
$\left[2\frac{\th_1}{\pi}\right]=2$. So
$\left[2\frac{\th_2}{\pi}\right]=3$ by (\ref{3.6}), which, together
with (\ref{3.7}), implies $\frac{3}{2}<\frac{\th_2}{\pi}<2$. And
hence $1<\frac{\th_1}{\pi}<\frac{7}{6}$ by (\ref{3.2}). In summary,
by the value of $i(c_2^4)$ we obtain the estimates
\be 1<\frac{\th_1}{\pi}<\frac{7}{6}\quad\mbox{and}\quad
\frac{3}{2}<\frac{\th_2}{\pi}<2.\lb{3.8}\ee

By Claim 1 we have
$4=i(c_2^5)=2\sum_{i=1}^{2}\left[\frac{5\th_i}{2\pi}\right]-8$,
which implies
\be \sum_{i=1}^{2}\left[\frac{5\th_i}{2\pi}\right]=6.\lb{3.9}\ee
Multiplying (\ref{3.8}) by $\frac{5}{2}$ yields
\be \frac{5}{2}<\frac{5\th_1}{2\pi}<\frac{35}{12}\quad\mbox{and}\quad
          \frac{15}{4}<\frac{5\th_2}{2\pi}<5. \lb{3.10}\ee
By (\ref{3.10}) we have $\left[\frac{5\th_1}{2\pi}\right]=2$. So
$\left[\frac{5\th_2}{2\pi}\right]=4$ by (\ref{3.9}), which, together
with (\ref{3.10}), implies $\frac{8}{5}<\frac{\th_2}{\pi}<2$. And
hence $1<\frac{\th_1}{\pi}<\frac{16}{15}$ by (\ref{3.2}). In
summary, by the value of $i(c_2^5)$ we obtain the estimates
\be 1<\frac{\th_1}{\pi}<\frac{16}{15}\quad\mbox{and}\quad
        \frac{8}{5}<\frac{\th_2}{\pi}<2. \lb{3.11}\ee

By Claim 1 we have
$4=i(c_2^6)=2\sum_{i=1}^{2}\left[3\frac{\th_i}{\pi}\right]-10$,
which implies
\be \sum_{i=1}^{2}\left[3\frac{\th_i}{\pi}\right]=7.\lb{3.12}\ee
Hence by (\ref{3.11}) we have
\be 3<3\frac{\th_1}{\pi}<\frac{16}{5}\quad\mbox{and}\quad
           \frac{24}{5}<3\frac{\th_2}{\pi}<6.\lb{3.13}\ee
By (\ref{3.13}) we have $\left[3\frac{\th_1}{\pi}\right]=3$. So
$\left[3\frac{\th_2}{\pi}\right]=4$ by (\ref{3.12}), which, together
with (\ref{3.13}), implies
$\frac{8}{5}<\frac{\th_2}{\pi}<\frac{5}{3}$. In summary, by the
value of $i(c_2^6)$ we obtain estimates
\be 1<\frac{\th_1}{\pi}<\frac{16}{15}\quad\mbox{and}\quad
         \frac{8}{5}<\frac{\th_2}{\pi}<\frac{5}{3}.\lb{3.14}\ee

By Claim 1 we have
$6=i(c_2^7)=2\sum_{i=1}^{2}\left[\frac{7\th_i}{2\pi}\right]-12$,
which implies
\be \sum_{i=1}^{2}\left[\frac{7\th_i}{2\pi}\right]=9.\lb{3.15}\ee
Hence by (\ref{3.14}) we have
\be \frac{7}{2}<\frac{7\th_1}{2\pi}<\frac{56}{15}<4\quad\mbox{and}\quad
       \frac{28}{5}<\frac{7\th_2}{2\pi}<\frac{35}{6}. \lb{3.16}\ee
By (\ref{3.16}) we have $\left[\frac{7\th_1}{2\pi}\right]=3$. So
$\left[\frac{7\th_2}{2\pi}\right]=6$ by (\ref{3.15}), which,
together with (\ref{3.16}), implies $6<\frac{7\th_2}{2\pi}<\frac{35}{6}$.
This leads to a contradiction.

{\bf Case  II}: {\it $c_1\in$NCG-2 with $i(c_1)=0$ or $1$.}

Then in this case $\hat{i}(c_1)$ is an irrational number by the
definition of $\th$. By the assumption, $c_1$ and $c_2$ do not
simultaneously belong to the classes
$\{\mbox{NCG-1},\mbox{NCG-2}\}$. Then $c_2\in\{\mbox{NCG-3,
NCG-4}\}$ with $\hat{i}(c_2)$ an integer, Hence
$\sum_{i=1}^2\frac{\ga_{c_i}}{\hat{i}(c_i)}$ is an irrational
number. However, by Proposition 2.3 we have
$\sum_{i=1}^2\frac{\ga_{c_i}}{\hat{i}(c_i)}=1$ contradiction!

\smallskip

{\bf Case III:} {\it $c_1\in$NCG-1 with $i(c_1)=0$.}

We continue our study 3 subcases:

(i) If $c_2\in\{\mbox{NCG-3, NCG-4}\}$ with $i(c_2)$ even, then the
condition of Lemma 3.3 is satisfied, i.e., $M_0=b_0=0$ by
Proposition 2.4. However, in this subcase we have $M_0\ge 1$,
contradiction!

(ii) If $c_2\in$NCG-3 with $i(c_2)=2p-1\ge3$, then in this subcase
we have $M_0\ge 1$ and $M_1=0$. But by Propositions 2.4 and 2.5 we
obtain $-1\ge M_1-M_0\ge b_1-b_0=0$, contradiction!

(iii) If $c_2\in$NCG-3 with $i(c_2)=1$, then this is exactly subcase
(v) in Case I.

This completes the proof of Lemma 3.4. \hfill\hb

{\bf Remark 3.5.} From the proof of Lemma 3.4, one can see that,
when there exist precisely two prime closed geodesics $c_1$ and $c_2$
on a bumpy $(S^3,F)$, at least one of them must have initial index $2$
provided they do not belong to the following two precise classes:

(1) $c_1\in$NCG-1 with $i(c_1)=0$ and $c_2\in$NCG-2 with $i(c_2)=1$.

(2) $c_1\in$NCG-2 with $i(c_1)=0$ and $c_2\in$NCG-2 with $i(c_2)=1$.

Note that Theorem 1.1 is a weaker consequence of the following
Theorem 3.6.

\smallskip

{\bf Theorem 3.6.} {\it Let $(M,F)$ be a bumpy Finsler $\Q$-homological $S^3$
with precisely two prime closed geodesics $c_1$ and $c_2$. Then both of $c_1$ and
$c_2$ must belong to classes $\{\mbox{NCG-1},\mbox{NCG-2}\}$. }

\smallskip

{\bf Proof of Theorem 3.6.} Because the proof is the same for $\Q$-homological
$S^3$, it suffices to prove Theorem 3.6 for bumpy Finsler $S^3=(S^3,F)$ only.
Suppose that there exist precisely two closed geodesics $c_1$ and $c_2$ on
$(S^3,F)$. Assume the theorem does not hold, without loss of generality, we can
suppose $i(c_1)=2$ by Lemma 3.4. Next we carry out our proof in three steps
according to the classification of $c_1$, and will derive some contradiction in
each case.

{\bf Step 1}: {\it $c_1\in\{\mbox{NCG-3, NCG-4}\}$}

In this case $i(c_1^m)=2m,\forall m\ge1$. Then by Proposition 2.1,
there holds $M_{2k}(1)=1,\,M_{2k-1}(1)=0,\forall k\in\N.$ Since
$M_{2k}=M_{2k}(1)+M_{2k}(2)\ge 2,\,\forall k\ge2$ by Propositions
2.4 and 2.5, it yields $M_{2k}(2)\ge 1,\,\forall k\ge2$. So $i(c_2)$
must be even by Proposition 2.1. We continue our study in 4
subcases:

(i) If $c_2\in\{\mbox{NCG-3},\;\mbox{NCG-4}\}$ with
$i(c_2^m)=2m,m\in\N$, then $M_2(2)=1$. Thus $2=M_2=b_2=1$ by
Proposition 2.4 and Lemma 3.3, contradiction!

(ii) If $c_2\in\{\mbox{NCG-3},\;\mbox{NCG-4}\}$ with $i(c_2)=2p\ge
4$, then we have $\ga_{c_1}=\ga_{c_2}=1$. Hence by Proposition 2.3
we obtain
$1=\sum_{k=1}^2\frac{\ga_(c_k)}{\hat{i}(c_k)}=\frac{1}{2}+\frac{1}{2p}<1$,
contradiction!

(iii) If $c_2\in\mbox{NCG-2}$, noting that $\hat{i}(c_2)$ is an
irrational number and $\hat{i}(c_1)$ is an integer, this leads to a
contradiction by Proposition 2.3.

(iv) If $c_2\in$NCG-1, then we have
$\hat{i}(c_2)=2(p-1)+\frac{\th_1+\th_2}{\pi}$. And by Definition
2.2, $\ga_{c_1}=\ga_{c_2}=1$. Hence by Proposition 2.3 we obtain
$\frac{1}{2(p-1)+\frac{\th_1+\th_2}{\pi}}+\frac{1}{2}=1$, i.e.,
$\frac{\th_1+\th_2}{\pi}=4-2p$. By the definitions of $\th_1$ and
$\th_2$, we have $\frac{\th_1+\th_2}{\pi}\in(0,4)$. So $p=1$, which
implies $M_2(2)\ge1$. Thus $2\le M_2(1)+M_2(2)=M_2=b_2=1$,
contradiction!

{\bf Step 2:} {\it $c_1\in$NCG-2 with $i(c_1)=2$.}

In this case $i(c_1^m) =m+2[\frac{m\th}{2\pi}]+1, \forall m\ge1$.
Hence by Proposition 2.1, $M_{2k}(1)=1$ for some $k\in\N$ and
$M_{2j-1}(1)=0$ for $j\in\N$. We continue our study in 8 subcases:

(i) If $c_2\in\{\mbox{NCG-3},\;\mbox{NCG-4}\}$, noting that
$\hat{i}(c_1)$ is an irrational number and $\hat{i}(c_2)$ is an
integer, this leads to a contradiction by Proposition 2.3.

(ii) If $c_2\in$NCG-2 with $i(c_2)$ odd, then $M_{2k}(2)=0, k\in\N$.
So we have $M_{2k}\le1,\forall k\ge2$, which contradicts to
$M_{2k}\ge b_{2k}=2,\forall k\ge2$ by Propositions 2.4 and 2.5.

(iii) If $c_2\in$NCG-2 with $i(c_2)=0$, then $M_0\ge1$ and
$M_{2k-1}=0, \forall k\in\N$. So by Proposition 2.4 and Lemma 3.3 we
obtain $1\le M_0=b_0=0$, contradiction!

(iv) If $c_2\in$NCG-2 with $i(c_2)=2$, then $M_2=2$ and $M_{2k-1}=0,
\forall k\in\N$. So by Proposition 2.4 and Lemma 3.3 we obtain $2=
M_2=b_2=1$, contradiction!

(v) If $c_2\in$NCG-2 with $i(c_2)=2p\ge 4$, $i(c_2^m) = m(2p-1)
+2[\frac{m\th}{2\pi}]+1, \forall m\ge1$. Let
$2T=2(2p-1)+2[\frac{\th}{\pi}]+2$. Then
$2<2T\notin\{i(c_2^m)\,|\,m\in\N\}$. So by Proposition 2.1
$M_{2T}(2)=0$, which implies $M_{2T}=M_{2T}(1)+M_{2T}(2)\le1$, where
$2T>2$. By Lemma 3.3 it yields $1\ge M_{2T}=b_{2T}=2$,
contradiction!

(vi) If $c_2\in$NCG-1 with $i(c_2)=0,i(c_2^m) = -2m
+2\sum_{i=1}^{2}\left[\frac{m\th_i}{2\pi}\right]+2,m\in\N$, then by
Proposition 2.1 and Lemma 3.3, we have $1\le M_0=b_0=0$,
contradiction!

(vii) If $c_2\in$NCG-1 with $i(c_2)=2$, then $M_2=2$ and
$M_{2k-1}=0, \forall k\in\N$. So by Proposition 2.4 and Lemma 3.3 we
obtain $2= M_2=b_2=1$, contradiction!

(viii) If $c_2\in$NCG-1 with $i(c_2)=2k\ge 4$, then $i(c_2^m) =
2m(k-1) +2\sum_{i=1}^{2}\left[\frac{m\th_i}{2\pi}\right]+2, \forall
m\ge1$ and $M_{2j-1}(2)=0,\forall j\in\N$. Let
$m_0=\min\{m\in\N|\;\sum_{i=1}^{2}\left[\frac{m\th_i}{2\pi}\right]\ge
1\}$. Because $\sum_{i=1}^{2}[\frac{\th_i}{2\pi}]=0$ by the
definitions of $\frac{\th_i}{\pi}$, $m_0\ge2$. Then
\bea i(c_1^{m_0-1})&=&2(m_0-1)(k-1)+2.\\
i(c_1^{m_0})&=&2m_0(k-1)+2+2\sum_{i=1}^{2}\left[\frac{m_0\th_i}{2\pi}\right]\nn\\
&\ge& i(c_1^{m_0-1})+4. \eea

\noindent Let $2T=i(c_1^{m_0-1})+2$. Notice that $i(c_1^m)$ is
non-decreasing, $2T\notin\{i(c_2^m),m\in\N\}$. So by Proposition 2.1
we have $M_{2T}(2)=0$, where $2T\ge 4$. So we have $M_{2T}\le1$. By
Lemma 3.3 we have $1\ge M_{2T}=b_{2T}=2$, contradiction!

{\bf Step 3:} {\it $c_1\in$NCG-1 with $i(c_1)=2$.}

In this case, we have $\hat{i}(c_1)=\frac{\th_1}{\pi}+\frac{\th_2}{\pi}$ and
\be i(c_1^m)=2\sum_{i=1}^{2}\left[\frac{m\th_i}{2\pi}\right]+2, \forall
       m\ge1\qquad\mbox{and}\qquad M_{2k-1}(1)=0,\forall k\in\N. \lb{3.19}\ee
We continue our study in 4 subcases:

(i) If $c_2\in\{\mbox{NCG-3,NCG-4}\}$ with $i(c_2^m)=2m,m\in\N$,
then $M_2=M_2(1)+M_2(2)=2$, which is a contradiction to Proposition
2.4 and Lemma 3.3.

(ii) If $c_2\in\{\mbox{NCG-3,NCG-4}\}$ with $i(c_2^m)=2pm\ge
4,m\in\N$ and $\hat{i}(c_2)=2p$, then we have
\be M_{2p}(2)=1,\quad M_q(2)=0,\;0\le q\le2p-1,\qquad M_{2k-1}(2)=0, \quad
      \forall k\in\N.\lb{3.20}\ee
By Definition 2.2 we have
$\ga_{c_1}=\ga_{c_2}=1$. So it follows from Proposition 2.3 in
Section 2 that
$\frac{1}{\frac{\th_1}{\pi}+\frac{\th_2}{\pi}}+\frac{1}{2p}=1$,
i.e., we have
\be \frac{2p-1}{2}(\frac{\th_1}{\pi}+\frac{\th_2}{\pi})=p.\lb{3.21}\ee
Noting that both $\frac{(2p-1)\th_1}{2\pi}$ and
$\frac{(2p-1)\th_2}{2\pi}$ are irrational numbers, we obtain
$[\frac{(2p-1)\th_1}{2\pi}]+[\frac{(2p-1)\th_2}{2\pi}]\in\{0,1,\cdots,p-1\}$.
Hence we have
\be i(c_1^{2p-1})=2\left(\left[\frac{(2p-1)\th_1}{2\pi}\right]
   +\left[\frac{(2p-1)\th_2}{2\pi}\right]\right)+2\in\{2,4,\cdots,2p\}.\lb{3.22}\ee

\smallskip

{\bf Claim 2}: $i(c_1^{2p-2})=i(c_1^{2p-1})=2p$.

In fact, by Proposition 2.4, Lemma 3.3 and (\ref{3.20}), we have
\be M_2(1)=1, \;\mbox{and if }2p-2\ge 4,\mbox{ then }
           M_{2q}(1)=2,\;4\le 2q\le2p-2. \lb{3.23}\ee

If $2p=4$, then by (\ref{3.19}) and (\ref{3.21}) it yields
$i(c_1^2)\le4$. We claim $i(c_1^2)=4$. Otherwise, assume
$i(c_1^2)=2$, then it yields $M_2=M_2(1)=2$, contradicting
(\ref{3.23}). So we have $i(c_1^2)=4$. Since $i(c_1^{3})\ge i(c_1^{2})=4$,
by (\ref{3.22}) we obtain $i(c_1^{3})=4$.

If $2p-2\ge4$, noting that $i(c_1^m)$ is non-decreasing, then by
(\ref{3.23}), $M_2$ is uniquely contributed by $c_1$, and
$M_{2q}(1)$ in (\ref{3.23}) should be uniquely contributed by the
two successive iterations $c_1^{m_0}$ and $c_1^{m_0+1}$ with
$i(c_1^{m_0})=i(c_1^{m_0+1})=2q$ for some $m_0\in\N$. So we have the
following sequence about the values of $i(c_1^m),1\le m\le 2p-2$
\bea &&i(c_1)=2,\\
&&i(c_1^{2j})=i(c_1^{2j+1})=2j+2,\;\mbox{for}\,j=1,\cdots,p-2,\\
&&i(c_1^{2p-2})=2p.\lb{3.26}\eea
Since $i(c_1^{2p-1})\ge i(c_1^{2p-2})$, by (\ref{3.22}) and (\ref{3.26}) we
obtain $i(c_1^{2p-1})=2p$. Claim 2 is proved.

\smallskip

By Claim 2, it yields $M_{2p}(1)\ge2$. So
$M_{2p}=M_{2p}(1)+M_{2p}(2)\ge3$, where $2p\ge4$. However, by
Proposition 2.4 and Lemma 3.3 we have $3=M_{2p}=b_{2p}=2$,
contradiction!

(iii) If $c_2\in\{\mbox{NCG-2},\;\mbox{NCG-3}\}$ with $i(c_2)$ odd,
then $M_{2k}(2)=0,k\in\N$ by Proposition 2.1. But by Propositions
2.4 and 2.5, we have $M_2\ge1$ and $M_{2k}\ge 2,\forall k\ge2$.
Hence Morse type numbers $M_{2k}$ must be contributed by iterations
of $c_1$, i.e., $M_{2k}=M_{2k}(1),k\in\N$. Thus $M_2$ should be
contributed at least by $c_1$, and $M_{2k},k\ge2$ should be
contributed at least by the two successive iterations $c_1^{m_0}$
and $c_1^{m_0+1}$ with $i(c_1^{m_0})=i(c_1^{m_0+1})=2k$ for some
$m_0\in\N$. Because $i(c_1^m)$ is non-decreasing,  we have
\be i(c_1^{k+2})-i(c_1^{k})\in\{0,2\},\;\forall k\in\N. \lb{3.27}\ee
But by the common index jump Theorem 3.1, there exist infinitely many
$k'\in\N$ such that
\be i(c_1^{2k'+1})-i(c_1^{2k'-1})=2i(c_1)=4, \lb{3.28}\ee
contradicting (\ref{3.27}).

(iv) If $i(c_2)\in$NCG-2 with $i(c_2)=2p\ge 4, i(c^m) = m(2p-1)
+2\left[\frac{m\th}{2\pi}\right]+1,m\in\N$ and
$\hat{i}(c_2)=2p-1+\frac{\th}{\pi}$. By Definition 2.2
$\ga_{c_1}=\ga_{c_2}=1$. So by Proposition 2.3, we have
$\frac{1}{\frac{\th_1}{\pi}+\frac{\th_2}{\pi}}+\frac{1}{2p-1+\frac{\th}{\pi}}=1$,
i.e.,
\be \frac{\th_1}{\pi}+\frac{\th_2}{\pi}<\frac{2p-1}{2p-2}\le\frac{3}{2}.\lb{3.29}\ee

By (\ref{3.29}) we obtain
$[\frac{\th_1}{\pi}]+[\frac{\th_2}{\pi}]\in\{0,1\}$. If
$[\frac{\th_1}{\pi}]+[\frac{\th_2}{\pi}]=0$, then
$i(c_1^2)=2([\frac{\th_1}{\pi}]+[\frac{\th_2}{\pi}]) +2=2$.  So
$M_2=M_2(1)=2$. However, by Proposition 2.4 and Lemma 3.3 we have
$2=M_2=b_2=1$, contradiction! Hence $[\frac{\th_1}{\pi}]+[\frac{\th_2}{\pi}]=1$,
which together with (\ref{3.29}) yields, without loss of generality,
\be 1<\frac{\th_1}{\pi}<\frac{2p-1}{2p-2}.\lb{3.30}\ee

 {\bf Claim 3}: {\it $i(c_1^{2p-1})=2p+2$ and $i(c_1^{4p-2})=4p$.}

In fact, in this subcase we have $M_{2k-1}=0,k\in\N$ by Proposition 2.1. So the
condition of Lemma 3.3 is satisfied. Noting that
$i(c_2)=2p, i(c_2^2)\in2\Z-1$ and $i(c_2^m)\ge 3(2p-1)+1>4p,\forall
m\ge3$, by Proposition 2.1 we obtain
\be M_{2p}(2)=1,\; M_q(2)=0,\; \forall\ 0\le q\le 4p,\ q\neq 2p. \lb{3.31}\ee
Hence Proposition 2.4, Lemma 3.3 and (\ref{3.31}) yield
\be M_2(1)=M_{2p}(1)=1,\;M_{2q}(1)=2,\;\;\forall\ 4\le 2q\le 4p,\ 2q\neq 2p.\lb{3.32}\ee
Noting that $i(c_1^m)$ is non-decreasing, so by (\ref{3.32}) we have
the following
\bea
&&\{i(c_1),i(c_1^2),i(c_1^3)\cdots,i(c_1^{x-3}),i(c_1^{x-2}),
i(c_1^{x-1}),i(c_1^{x}),i(c_1^{x+1})\cdots i(c_1^{y-1}),i(c_1^{y})\}\nn\\
&&\;\;=\{2,4,4,\cdots,2p-2,2p-2,2p,2p+2,2p+2, \cdots,4p,4p\}, \lb{3.33}\eea
where $x=\min\{m\in\N\,|\,i(c_1^m)=2p+2\}$ and $y=\max\{m\in\N\,|\,i(c_1^m)=4p\}$,
which are determined by the equations
\bea
x&=&2\cdot\frac{2p-4}{2}+3=2p-1.\lb{3.34}\\
y&=&2\cdot\frac{2p-4}{2}+1+2\cdot\frac{4p-2p}{2}+1=4p-2.\lb{3.35}\eea
So by (\ref{3.33}), (\ref{3.34}) and (\ref{3.35}), Claim 3 is
proved.

\smallskip

By (\ref{3.30}), we obtain $\frac{4p-2}{2}\frac{\th_1}{\pi}=(2p-1) \frac{\th_1}{\pi}>2p-1$.
So by Claim 3 and (\ref{3.19}), it yields
\be 4p=i(c_1^{4p-2})=2\sum_{i=1}^{2}\left[\frac{(4p-2)\th_i}{2\pi}\right]+2\ge4p. \lb{3.36}\ee
Thus $\frac{4p-2}{2}\frac{\th_1}{\pi}<2p$, i.e.,
\be \frac{\th_1}{\pi}<\frac{2p}{2p-1}.\lb{3.37}\ee

On the other hand, by (\ref{3.29}) and (\ref{3.30}), it yields
$\frac{\th_2}{\pi}<\frac{2p-1}{2p-2}-1=\frac{1}{2p-2}$. Since
$2p\ge4$, we have \be
\frac{(2p-1)\th_2}{2\pi}<\frac{2p-1}{2(2p-2)}<1.  \lb{3.38}\ee

\noindent Thus by Claim 3 and (\ref{3.19}), it yields
\be 2p+2=i(c_1^{2p-1}) =2\sum_{i=1}^{2}\left[\frac{(2p-1)\th_i}{2\pi}\right]+2
         =2\left[\frac{(2p-1)\th_1}{2\pi}\right]+2,  \lb{3.39}\ee
which implies $p<\frac{2p-1}{2}\frac{\th_1}{\pi}<p+1$. This is to say
\be \frac{\th_1}{\pi}>\frac{2p}{2p-1}. \lb{3.40}\ee
contradicting (\ref{3.37}).

Above three steps complete the proof of Theorem 3.6. \hfill\hb

{\bf Remark 3.7.} Suppose that there exist precisely two prime closed geodesics
$c_1$ and $c_2$ on a bumpy $(S^3,F)$. From Remark 3.5 and the proof of
Theorem 3.6, both $c_1$ and $c_2$ are non-hyperbolic and must belong to
one of the following precise classes:

(I) $c_1\in$NCG-1 with $i(c_1)=0$ and $c_2\in$NCG-2 with $i(c_2)=1$.

(II) $c_1\in$NCG-2 with $i(c_1)=0$ and $c_2\in$NCG-2 with
$i(c_2)=1$.

(III) $c_1\in$NCG-1 with $i(c_1)=2$ and $c_2\in$NCG-1 with
$i(c_2)=2p\ge 4$.

\bibliographystyle{abbrv}

\end{document}